\documentclass[12pt]{article}

\usepackage[small,nohug,heads=vee]{diagrams}
\diagramstyle[labelstyle=\scriptstyle]

\usepackage[latin1]{inputenc}
\usepackage[T1]{fontenc}
\usepackage{hyperref} 
\usepackage[all]{xy}
\usepackage{amsfonts}
\usepackage{amsmath,amssymb}
\usepackage{amscd}
\usepackage{amsthm}
\usepackage{slashed}
\newtheorem{df}{Definition}[section]
\newtheorem{thm}{Theorem}[section]
\newtheorem{lem}{Lemma}[section]
\newtheorem{cly}{Corollary}[section]
\newtheorem{prop}{Proposition}[section]
\newtheorem{rem}{Remark}[section]

\newtheorem*{mycly}{Corollary}

\newcommand{\g}{\mathfrak{g}}
\newcommand{\h}{\mathfrak{h}}
\newcommand{\e}{\mathfrak{e}}
\newcommand{\el}{\mathfrak{l}}

\title{Non abelian cohomology of extensions of Lie algebras as Deligne groupoid}
\author {Ya\"el Fr\'egier \\
MIT/Z\"urich/Lens }
\begin {document}
\date{}
\maketitle

\begin{abstract}
In this note we show that the theory of non abelian extensions of a Lie algebra $\g$ by a   Lie algebra $\h$ can be understood in terms of a differential graded Lie algebra $L$. More precisely we show that the non-abelian cohomology $H^2_{nab}(\g,\h)$ is  $\mathcal{MC}(L)$, the $\pi_0$ of the Deligne groupoid of $L$.
\end{abstract}

\tableofcontents

\section*{Introduction}
\addcontentsline{toc}{section}{Introduction}

The aim of this note is to show that equivalence classes of extensions of a Lie algebra $\g$ by a Lie algebra $\h$ can be seen as the $\pi_0$ of the Deligne groupoid of a {\bf d}ifferential {\bf g}raded {\bf L}ie {\bf a}lgebra ({\bf dgLa}) $(L, [\ ,\ ])$.
In other words, Theorem \ref{princip} states that
$$H^2_{nab}(\g,\h)\simeq \mathcal{MC}(L).$$

Our first motivation is that the  formalism of {\bf dgLa} enables to extend notions from classical algebraic structures to their up to homotopy and geometrical counterparts. 

More precisely, Lie algebras on a vector space $V$, can be seen as elements $\alpha$ of degree one in some graded Lie algebra $(L, [\ , \, ])$ satisfying $[\alpha,\alpha]=0$ (\cite{NR}).  Here $(L, [\ ,\ ])$ depends only on V.  The same equation in a setting where V is replaced by a graded vector space  leads to the notion of homotopy Lie algebra (\cite{LS}).  

Most of the algebraic notions related to Lie algebras like modules, semi-direct products, cohomology, etc. can be defined for homotopy Lie algebras. A way to do it is to describe the Lie algebra version of these notions in terms of {\bf dgLa}. Considering the same formulation in a graded setting gives their homotopy version.

In particular, the aim of our work \cite{BF} is to build a theory of extensions of homotopy Lie algebras. This is why we first need to understand the theory of extensions of Lie algebras in terms of {\bf dgLa}.

Analogously, geometrical objects generalizing Lie algebras like Lie algebroids, Poisson tensors, Courant algebroids, etc. admit a similar description in terms of {\bf dgLa}. So the present work should open the door to a theory of non-abelian extensions in these categories.\\

Our second motivation is the absence of a systematic treatment of non abelian cohomology of extensions of algebras. The correct setting to give such a theory should be the language of operads. This setting, thanks to Koszul duality  (\cite{GK}), furnishes to any Koszul operad  graded Lie algebras similar to $(L,[ \ ,\ ])$. Therefore, understanding non-abelian cohomology and extensions in terms of  {\bf dgLa} is the key step leading to a systematic theory via the operadic  machinery. \\

In this paper we start by recalling section \ref{background} the classical theory of classification of extensions of algebras in terms of non-abelian cohomology. We also give the definitions of Deligne groupoid necessary to formulate our main result Thm. \ref{princip}.

 We then turn in section \ref{ArtinDeligne} to the proof of this result. We recall the Chevalley-Eilenberg cohomology and the Nijenhuis-Richardson bracket in order to define the {\bf dgLa} $(L,[\ ,\ ])$ (Def \ref{CE}, Lemma \ref{L}). We show in section \ref{MC} that extensions are in bijections with Maurer-Cartan elements. We show in section \ref{gauge}  (Prop. \ref{equiv}) that the equivalence relation for non-abelian cohomology coincides with gauge equivalence for Maurer-Cartan elements.

Finally we give in section \ref{tangent} (Prop. \ref{tangentcomplex}) the relation, in terms of tangent cohomology complex, between A) the classification of abelian extensions of a Lie algebra in terms of Chevalley-Eilenberg cohomology and B) the Deligne groupoid.

\newpage

\section{Theoretical background}\label{background}

\subsection{Deligne groupoid}

The purpose of this section is to recall the classical notion of Deligne groupoid of a differential graded Lie algebra, as can be found for example in \cite{GM},  or in \cite{LV} 13.3.20. We work over a field $\mathbb{K}$ of characteristic zero.
\begin{df} A \emph{graded Lie algebra} is a $\mathbb{Z}$-graded vector space $L=\bigoplus_{n\in\mathbb{Z}}L_n$ equipped with a degree-preserving bilinear bracket $[
\cdot,\cdot] \colon L\otimes L \longrightarrow L$ which satisfies 
\begin{itemize}
\item[1)] graded antisymmetry: $[a,b]=-(-1)^{\vert a\vert\vert b\vert}[b,a]$,
\item[2)] graded Leibniz rule: $[a,[b,c]]=[[a,b],c] + (-1)^{\vert a\vert\vert b\vert}[b,[a,c]].$
\end{itemize}
Here $a, b, c$ are homogeneous elements of $L$ and the degree $\vert x\vert$ of an homogeneous element $x\in L_n$ is by definition $n$. 
\end{df}

 \begin{df} A \emph{differential graded Lie algebra} ({\bf dgLa} for short) is {a} graded Lie algebra $(L,[\cdot,\cdot])$ equipped with a homological derivation $d \colon L\to L$ of degree 1. In other words:
\begin{itemize}
\item[1)] $\vert da\vert=\vert a\vert +1$ ($d$ of degree 1),  
\item[2)] $d[a,b]=[da,b]+ (-1)^{ \vert a\vert}[a,db]$ (derivation),
\item[3)] $d^2=0$ (homological).
\end{itemize}
\end{df}

\begin{df} The set $MC(L)$ of Maurer-Cartan elements of a {\bf dgLa} $(L,[\cdot,\cdot],d)$ is defined to be 
$$\{\alpha \in L^1\mid d\alpha+\frac{1}{2}[\alpha,\alpha]=0\}.$$

\end{df}

We now want to recall the equivalence relation on $MC(L)$, coming from {\it gauge symmetries}. We motivate the formulas by embedding the {\bf dgLa} $L$ into a graded Lie algebra $L_s$.

\begin{df}
Let $(L,[\ ,\ ],d)$ be a differential graded Lie algebra. A specter of L consists in a graded Lie algebra $(L_s, [\ ,\ ]_s)$ containing L as a sub-Lie algebra together with a degree one element $D$ such that $[D,D]_s=0$ and $d=[D, \ ]_s$.
\end{df}

One can remark that since $d=[D, \ ]_s$, $D$ lies in the normalizer $N(L)$ of $L$.\\

The symmetries of $L$ that we want to consider are simply the shadows of the gauges symmetries of $L_s$ : the action of $\beta \in L^0$ on the element $l_s := l+D$  reads 
$$e^\beta l_s=e^\beta l+e^\beta D=e^\beta l+(e^\beta-1 )D +D,$$
where ``$e^\beta$'' means $e^{ad_\beta }=e^{[\beta ,\ ]}$
\begin{rem}
We assume here and in the following that $\beta$ is ad-nilpotent, i.e. $\forall x\in L, \exists n/ {ad_\beta}^n(x)=0$. This in order to have a convergent sum defining $e^\beta$. It  will be the case in the applications we consider.
\end{rem}
Since $D$  is in $N(L)$, the term $$g_\beta^s:=(e^\beta-1 )D$$ lies in $L$.

We want to have an intrinsic definition of $g_\beta^s$ , i.e. a direct expression without invoking any specter of $L$. Let us therefore define the function $\slashed{e}$ by $$\slashed{e}^x:=\sum_{n\in\mathbb{N}^*}\frac{1}{(n+1)!}x^n.$$

\begin{df}
$(L,[\ ,\ ],d)$ be a differential graded Lie algebra, and let $\beta$ in $L^0$ be ad-nilpotent. One defines
\begin{equation}
g_\beta:=-\slashed{e}^{ad_\beta} d\beta.\label{gbeta}
\end{equation}
\end{df}

\begin{lem}
$$g_\beta^s=g_\beta.$$
\end{lem}

\begin{proof}
\begin{eqnarray*}
g_\beta^s & = & \sum_{n\in \mathbb{N^*}} \frac{1}{n!}{ad_\beta}^nD\\
                   & = & \sum_{n\in \mathbb{N^*}} \frac{1}{n!}{ad_\beta}^{n-1}[\beta,D]_s\\
                   & = & - \sum_{n\in \mathbb{N}} \frac{1}{(n+1)!}{ad_\beta}^{n}d\beta\\
                   & = & g_\beta
                   \end{eqnarray*}
\end{proof}

We can introduce the relation $\sim$:

\begin{df}\label{equiMC}
$(L,[\ ,\ ],d)$ be a differential graded Lie algebra, and let $l$ and $l'$ be two elements in $L$. One says that $l$ is in relation with $l'$ and write $$l\overset{\beta}{\sim} l'$$ to signify that there exists an ad-nilpotent $\beta$  in $L^0$ such that
\begin{equation}\label{action}
l'=e^\beta l + g_\beta. 
\end{equation}
\end{df}

\begin{lem}
If $L^0$ is abelian, the relation $\sim$ defines an equivalence relation on $L^1$. 
\end{lem}

\begin{proof} 
 It is proven in section 1.3  of \cite{GM}, pages 50-51 that an arbitrary simply connected group $\mathcal{L}$ with Lie algebra $L_0$ acts on $L^1$ via the formula (\ref{action}). 
\end{proof}

\begin{lem}\label{stab}
The relation $\sim$ induces an equivalence relation on $MC(L)$.
\end{lem}

\begin{proof}
Once again we refer to \cite{GM} (last paragraph of section 1.3).
\end{proof}

\begin{df}
$(L,[\ ,\ ],d)$ be a differential graded Lie algebra, with $L^0$ abelian. One defines the {\it  Deligne groupoid} of L to be the groupoid whose set of objects is $MC(L)$ and whose set of homs are empty for non equivalent objects and reduced to one element for equivalent objects. One defines its set of connected components $\pi_0$:
 
$$\mathcal{MC}(L):=\pi_0(MC(L)):=MC(L)/\sim.$$ 
\end{df}

\begin{rem}
Note that the notion of Deligne groupoid makes sense for more general differential graded Lie algebras. However the price to pay is to use the concept of  deformation functor on Artin rings. Our problem does not require this generality.
\end{rem}

%%%%%%%%%%%%%%%%%%%%%%%%%%%%%%%%%%%%%%%%%%%%%%%%

\subsection{Extensions of Lie algebras classified by non abelian cohomology}\label{artin}
%%%%%%%%%%%%%%%%%%%%%%%%%%%%%%%%%%%%%%%%%%%%%%%%

Classically, extensions of algebras are defined in terms of short exact sequences.

\begin{df}
Let $\g$ and $\h$ be two Lie algebras. An extension $\e$ of $\g$ by $\h$ is a short exact sequence of the form
$$0\to \h\to \e \to \g \to 0.$$
\end{df}

But one is really interested in classes of extensions modulo the following equivalence relation :

\begin{df}
Let $\e$ and $\e'$ be two extensions of $\g$ by $\h$. They are said to be equivalent if there exists a commutative diagram 

  $$\xymatrix{
      0 \ar[r]  & \h \ar[r] \ar@{=}[d]& \e \ar[r] \ar[d]_\varphi&  \g  \ar[r] \ar@{=}[d]& 0  \\
   0 \ar[r] & \h \ar[r] & {\e^{'}} \ar[r] &  \g  \ar[r] & 0. }$$

\end{df}

Classes of extensions of $\g$ by $\h$ are classified in terms of nonabelian cohomology.

\begin{df}\label{nab}
A nonabelian 2-cocycle on $\g$ with values in $\h$ is a couple $(\chi, \psi)$ of linear maps $\chi : \g\wedge \g\to \h$ and $\psi : \g\to Der(\h)$ satisfying
\begin{equation}\label{mod}
[\psi(a),\psi(b)]=\psi([a,b])+ad_\h(\chi(a\wedge b))
\end{equation}
and
\begin{equation}\label{cocy}
\sum_{\circlearrowleft}\psi_a\chi(b\wedge c)-\chi([a,b]\wedge c)=0,
\end{equation}
where the sum is taken over cyclic permutations of $a$,$b$ and $c$.
One denotes by $Z^2_{nab}(\g,\h)$ this set of nonabelian cocycles. 

Nonabelian cohomology $H^2_{nab}(\g,\h)$ will be the quotient of $Z^2_{nab}(\g,\h)$ by the equivalence relation : $(\chi, \psi)\overset{\beta}{\thicksim}(\chi', \psi')$  if there exists $\beta: \g \to \h$ satisfying
\begin{equation}\label{eq1}
\psi'_a=\psi_a+ad^\h_{\beta(a)}
\end{equation}
and 
\begin{equation}\label{eq2}
\chi'(a\wedge b)=\chi(a\wedge b)+\psi_a\beta(b)-\psi_b\beta(a)-\beta([a,b])+[\beta(a),\beta(b)].
\end{equation}

\end{df}

\begin{prop}\label{class}
Extensions of $\g$ by $\h$ are classified by $H^2_{nab}(\g,\h)$.
\end{prop}

\begin{proof}
One associates a class in $H^2_{nab}(\g,\h)$ to an extension $\e$ through the choice of a section, i.e. a map $s : \g\to \e$ in
$$  \xymatrix{ 0 \ar[r] & \h \ar[r] & \e \ar[r]_p &  \g \ar@/_/@{.>}[l]_{s}  \ar[r] & 0},$$
  such that $p\circ s=Id_\g$.

 The representative cocycle $(\psi^s, \chi^s)$ is defined by the formulas
\begin{equation}
\psi^s_g(h):=[s(g),h]_\e
\end{equation}
and
\begin{equation}
\chi^s(g_1,g_2):=[s(g_1),s(g_2)]_\e-s([g_1,g_2]_\g).
\end{equation}

The point of the equivalence relation on $Z^2_{nab}(\g,\h)$ is to make this assignment independent of the choice of the section. Two sections $s_1$ and $s_2$ differ by an element $\beta : \g\to \h$ and one can check that $(\psi^{s_1}, \chi^{s_1})\overset{\beta}{\thicksim}(\psi^{s_2}, \chi^{s_2})$, with $\beta:=s_1-s_2$.

Moreover, this assignment sends equivalent extensions to the same element in $H^2_{nab}(\g,\h)$. Consider the following equivalence of extensions, with chosen sections $s$ and $s'$

  $$\xymatrix{
      0 \ar[r]  & \h \ar[r] \ar@{=}[d]& \e \ar[r] \ar[d]_\varphi&  \g  \ar@/_/@{.>}[l]_{s}\ar[r] \ar@{=}[d]& 0  \\
   0 \ar[r] & \h \ar[r] & {\e^{'}} \ar[r] &  \g \ar@/ ^/@{.>}[l]^{s'} \ar[r] & 0,}$$
  
then one has $(\psi^{s}, \chi^{s})\overset{\beta_\varphi}{\thicksim}(\psi^{s'}, \chi^{s'})$, with $\beta_\varphi:=\varphi^{-1}\circ s'-s$.\\

 In the other direction, on can associate to a cocycle $(\psi, \chi)$ an extension of the form 
 $$0\to \h\to \g\oplus\h_{(\psi, \chi)} \to \g \to 0, $$
 with bracket on $\g\oplus\h$ defined by
$$[h_1+g_1,h_2+g_2]_\e=[h_1,h_2]_\h+\psi_{g_1}(h_2)-\psi_{g_2}(h_1)+\chi(g_1,g_2)+[g_1,g_2]_\g.$$

This assignment is well defined at the level of cohomology classes since two equivalent cocycles $(\chi, \psi)\overset{\beta}{\thicksim}(\chi', \psi')$ lead to equivalent extensions

 $$\xymatrix{
      0 \ar[r]  & \h \ar[r] \ar@{=}[d]& \g\oplus\h_{(\psi, \chi)} \ar[r] \ar[d]_\varphi&  \g  \ar[r] \ar@{=}[d]& 0  \\
   0 \ar[r] & \h \ar[r] & {\g\oplus\h_{(\psi', \chi')}} \ar[r] &  \g  \ar[r] & 0,}$$

\end{proof}

More details can be found in \cite{AMR} and \cite{Nlab}.

%%%%%%%%%%%%%%%%%%%%%%%%%%%%%%%%%%%%%%%%%%%%%%%%%%%%%%%%%%%%%%%%%%%%%%%%%%%%%%%%%%%%%%%%%%%%%%%%%%%%

\newpage

%%%%%%%%%%%%%%%%%%%%%%%%%%%%%%%%%%%%%%%%%%%%%%%%%%%%%%%%%%%%%%%%%%%%%%%%%%%%%%%%%%%%%%%%%%%%%%%%%%%%

\section{Extensions in terms of Deligne groupoids} \label{ArtinDeligne}

The main result of this paper is the interpretation of $H^2_{nab}(\g,\h)$, the non-abelian 2nd cohomology of $\g$ with values in $\h$, as the Deligne groupoid of a suitable {\bf dgLa} $L$. 

Therefore we introduce in section  \ref{MC}  the {\bf dgLa} $L$ (lemma \ref{L}) and show that the set of its Maurer-Cartan elements is in bijection with the set of extensions of $\g$ by $\h$. In particular, thanks to the analysis done in the previous section, one obtains the  corollary \ref{Z2} :

\begin{mycly}
 $$Z^2_{nab}(\g,\h)\simeq MC(L).$$
\end{mycly}

We then show in section \ref{gauge} that the notion of gauge equivalence coincides with the notion of equivalence of non-abelian cohomology. Hence one obtains our main theorem, as a conjunction of corollary \ref{Z2} and of proposition \ref{equiv} :

\begin{thm} \label{princip} With the notations of definition  \ref{nab}
$$H^2_{nab}(\g,\h)\simeq \mathcal{MC}(L).$$
\end{thm}

%%%%%%%%%%%%%%%%%%%%%%%%%%%%%%%%%%%%%%%%%%%%%%%%%%%%%%%%%%%%%%%%%%%%%%%%%%%%%%%%%%%%%%%%%%%%%%%%%%%%

\subsection{Non abelian cocycles as Maurer-Cartan elements}\label{MC}

We start by analyzing the admissible forms that can have the Lie  bracket $\rho$ of an extension
$\e$ of $\g$ by $\h$

$$\xymatrix{
      0 \ar[r]  & \h \ar[r] \ar@{=}[d]& \e \ar[r]_{pr} \ar@{=}[d] & \g  \ar@/_/@{.>}[l]_{s}\ar[r] \ar[d]_\wr& 0  \\
   0 \ar[r] & H \ar[r] & {G\oplus H} \ar[r] &  G  \ar[r] & 0,}$$
 where $H$ is the vector space image of $\h$ in $\e$ and $G=s(\g)$ an arbitrary supplementary in $\e$.

 We first need to introduce a handy notation. Consider the canonical projection $$P_{A_1\dots A_n} : (G\oplus H)^{\otimes^n}\longrightarrow A_1\otimes \dots \otimes A_n,$$ where $A_i\in \{H,G\}$.
 
Let L be a linear map $L : (G\oplus H)^{\otimes^n}\longrightarrow   G\oplus H$ and let us denote
$$L_{A_1\dots A_n}^{A_{n+1}}:=P_{A_{n+1}}\circ L\circ P_{A_1\dots A_n}.$$
 With this notation and the ones of definition \ref{nab}, one has :
 \begin{rem} \label{form}

 All the components of $\rho$ are

$$\begin{array}{ c c c }
\begin{tabular}{ | c | c | c | }
\hline
  $\rho_{GG}^G$ &  $\rho_{GH}^G$ & $\rho_{HH}^G$   \\ \hline
 $\rho_{GG}^H$ &  $\rho_{GH}^H$  &  $\rho_{HH}^H$  \\ \hline
 \end{tabular} & = & \begin{tabular}{ | c | c | c | }
\hline
  $[\ ,\ ]_\g$ &   $0$ & $0$  \\ \hline
$\chi$ &  $\psi$   &  $[\ ,\ ]_\h$  \\ \hline
 \end{tabular}  \\
\end{array},$$

where $\rho_{GH}^G$ and $\rho_{HH}^G$ vanish since $P$ is a map of Lie algebras with kernel $H$, and one has identified $\rho_{GG}^G$ with $[\ ,\ ]_\g$ by conjugation by  $pr_{\vert G}$ whose inverse is $s$.
\end{rem}

Let us now analyze the different components of the jacobiator of $\rho$.

\begin{prop}  \label{Jac} The vanishing of the components of the Jacobiator of $\rho$ gives the following constraints :
\begin{enumerate}
\item[$J^G_{GGG}$] : $[\ ,\ ]_\g$ satisfies Jacobi identity,\\
\item[$J^H_{GGG}$] : $\chi$ is a ``Chevalley-Eilenberg'' cocycle (def \ref{CE}),\\
\item[$J^H_{GGH}$] : $H$ is a twisted $\g$-module,\\
\item[$J^H_{GHH}$] : $\g$ acts by derivations of $[\ ,\ ]_\h$,\\
\item[$J^H_{HHH}$] : $[\ ,\ ]_\h$ satisfies Jacobi identity.\\
\end{enumerate}
\end{prop}

The constraint coming from the vanishing of $J^H_{GGH}$, 
$$ [\psi(g_1),\psi(g_2)] =ad_\h(\chi(g_1\wedge g_2))+\psi([g_1,g_2]),$$
is exactly equation (\ref{mod}).
 It can be reformulated, if one denotes  by $\cdot$ the action defined by $\rho_{GH}^H,$ as
 $$g_1\cdot(g_2\cdot h)-g_2\cdot(g_1\cdot h)=ad_\h(\chi(g_1\wedge g_2)) +[g_1,g_2]_\g\cdot h.$$
 In particular, one recognizes the relation defining H as a $\g$-module, but twisted by the term $ad_\h(\chi(g_1\wedge g_2)).$\\

The use of the term Chevalley-Eilenberg ``cocycle'' is licit only when H is actually  $\g$-module, i.e. when the twisting $ad_\h(\chi(g_1\wedge g_2))$ vanishes.

\begin{proof} In this proof $t:=(e_1,e_2,e_3),$ with $e_i:=g_i+h_i\in \g\oplus\h.$
We have 
\begin{eqnarray*}
J^G_{GGG}(t) & =  & \sum_\circlearrowright (\rho^G_{GG}(e_{\sigma(1)},\rho^G_{GG}(e_{\sigma(2)},e_{\sigma(3)}))+\underset{0}{\underbrace{\rho^G_{GH}}}(e_{\sigma(1)},\rho^H_{GG}(e_{\sigma(2)},e_{\sigma(3)})))\\
& =  & \sum_\circlearrowright [g_1,[g_2,g_3]_\g]_\g.
\end{eqnarray*}

\begin{eqnarray*}
J^G_{GGH}(t) & =  & \sum_\circlearrowright (\rho^G_{GG}(e_{\sigma(1)},\underset{0}{\underbrace{\rho^G_{GH}}}(e_{\sigma(2)},e_{\sigma(3)})) +\underset{0}{\underbrace{\rho^G_{HG}}}((e_{\sigma(1)},\rho^G_{GG}(e_{\sigma(2)},e_{\sigma(3)}))  \\
 & +  &  \underset{0}{\underbrace{\rho^G_{HH}}}(e_{\sigma(1)},\rho^H_{GG}(e_{\sigma(2)},e_{\sigma(3)}))+ \underset{0}{\underbrace{\rho^G_{GH}}}(e_{\sigma(1)},\rho^H_{GH}(e_{\sigma(2)},e_{\sigma(3)})) )\\
& =  & 0.
\end{eqnarray*}

\begin{eqnarray*}
J^G_{GHH}(t) & =  & \sum_\circlearrowright (\rho^G_{GG}(e_{\sigma(1)},\underset{0}{\underbrace{\rho^G_{HH}}}(e_{\sigma(2)},e_{\sigma(3)}))+ \underset{0}{\underbrace{\rho^G_{HG}}}(e_{\sigma(1)},\rho^G_{GH}(e_{\sigma(2)},e_{\sigma(3)})) \\
& + &  \underset{0}{\underbrace{\rho^G_{GH}}}(e_{\sigma(1)},\rho^H_{HH}(e_{\sigma(2)},e_{\sigma(3)})) + \underset{0}{\underbrace{\rho^G_{HH}}}(e_{\sigma(1)},\rho^H_{GH}(e_{\sigma(2)},e_{\sigma(3)})))\\
& =  & 0.
\end{eqnarray*}

\begin{eqnarray*}
J^G_{HHH}(t) & =  &\sum_\circlearrowright ( \rho^G_{HG}(e_{\sigma(1)},\underset{0}{\underbrace{\rho^G_{HH}}}(e_{\sigma(2)},e_{\sigma(3)})) \\
& + & \underset{0}{\underbrace{\rho^G_{HH}}}(e_{\sigma(1)},\rho^H_{HH}(e_{\sigma(2)},e_{\sigma(3)})))\\
& =  & 0.
\end{eqnarray*}

\begin{eqnarray*}
J^H_{GGG}(t) & =  & \sum_{\sigma \in\circlearrowright}( \underset{\heartsuit_1}{\underbrace{\rho^H_{GG}(e_{\sigma(1)},\rho^G_{GG}(e_{\sigma(2)},e_{\sigma(3)}))}} + \underset{\heartsuit_2}{\underbrace{\rho^H_{GH}(e_{\sigma(1)},\rho^H_{GG}(e_{\sigma(2)},e_{\sigma(3)}))}})\\
& =  & \sum_{\sigma \in\circlearrowright}  (\chi (g_{\sigma(1)},[g_{\sigma(2)},g_{\sigma(3)}]_\g)+ g_{\sigma(1)}\cdot \chi (g_{\sigma(2)},g_{\sigma(3)}))\\
& =  & \delta \chi (g_1, g_2, g_3).
\end{eqnarray*}

\begin{eqnarray*}
J^H_{GGH}(t) & =  &\sum_{\sigma \in\circlearrowright} ( \rho^H_{GG}(e_{\sigma(1)},(\underset{0}{\underbrace{\rho^G_{GH}}}+\underset{0}{\underbrace{\rho^G_{HG}}})(e_{\sigma(2)},e_{\sigma(3)}))   \\
& +  &\underset{\Delta_1}{\underbrace{\rho^H_{HG}(e_{\sigma(1)},\rho^G_{GG}(e_{\sigma(2)},e_{\sigma(3)}))}}   +  \underset{\Delta_3^{''}}{\underbrace{\rho^H_{GH}(e_{\sigma(1)},\rho^H_{GH}(e_{\sigma(2)},e_{\sigma(3)}))}} \\
& +  &  \underset{\Delta_3^{'}}{\underbrace{\rho^H_{GH}(e_{\sigma(1)},\rho^H_{HG}(e_{\sigma(2)},e_{\sigma(3)}))}}   + \underset{\Delta_2}{\underbrace{\rho^H_{HH}(e_{\sigma(1)},\rho^H_{GG}(e_{\sigma(2)},e_{\sigma(3)}))}})\\
& =  & \sum_{\sigma \in\circlearrowright}(\underset{\Delta_3}{\underbrace{[\psi(g_{\sigma(1)}),\psi(g_{\sigma(2)})] }} -\underset{\Delta_2}{\underbrace{ ad_\h(\chi(g_{\sigma(1)}\wedge g_{\sigma(2)}))}})\\
& - & \underset{\Delta_1}{\underbrace{\psi([g_{\sigma(1)},g_{\sigma(2)}])}})(h_{\sigma(3)}) .
\end{eqnarray*}

\begin{eqnarray*}
J^H_{GHH}(t) & =  & \sum_{\sigma \in\circlearrowright} (\rho^H_{HG}(e_{\sigma(1)},(\underset{0}{\underbrace{\rho^G_{GH}}}+\underset{0}{\underbrace{\rho^G_{HG}}})(e_{\sigma(2)},e_{\sigma(3)}))    \\
& +  &\rho^H_{GG}(e_{\sigma(1)},\underset{0}{\underbrace{\rho^G_{HH}}}(e_{\sigma(2)},e_{\sigma(3)})) + \underset{g_{\sigma(1)}\cdot [h_{\sigma(2)},h_{\sigma(3)}]_\h=\Box_1}{\underbrace{\rho^H_{GH}(e_{\sigma(1)},\rho^H_{HH}(e_{\sigma(2)},e_{\sigma(3)}))} } \\
& +  & \underset{[-h_{\sigma(1)}, g_{\sigma(3)}\cdot h_{\sigma(2)}]_\h=\Box_3}{\underbrace{\rho^H_{HH}(e_{\sigma(1)},\rho^H_{HG}(e_{\sigma(2)},e_{\sigma(3)}))}}+\underset{-[g_{\sigma(2)}\cdot h_{\sigma(3)},h_{\sigma(1)}]_\h=\Box_2}{\underbrace{ \rho^H_{HH}(e_{\sigma(1)},\rho^H_{GH}(e_{\sigma(2)},e_{\sigma(3)})))}}\\
& =  &  \sum_{\sigma \in\circlearrowright} (g_{\sigma(1)}\cdot [h_{\sigma(2)},h_{\sigma(3)}]_\h-[g_{\sigma(1)}\cdot h_{\sigma(2)},h_{\sigma(3)}]_\h\\
& - & [h_{\sigma(2)}, g_{\sigma(1)}\cdot h_{\sigma(3)}]_\h).
\end{eqnarray*}

\begin{eqnarray*}
J^H_{HHH}(t) & =  & \sum_\circlearrowright \rho^H_{HG}(e_{\sigma(1)},\underset{0}{\underbrace{\rho^G_{HH}}}(e_{\sigma(2)},e_{\sigma(3)})) + \rho^H_{HH}(e_{\sigma(1)},\rho^H_{HH}(e_{\sigma(2)},e_{\sigma(3)}))\\
& =  & \sum_\circlearrowright [h_{\sigma(1)},[h_{\sigma(2)},h_{\sigma(3)}]_\h]_\h.
\end{eqnarray*}
\end{proof}

One recalls the definition of the Chevalley-Eilenberg complex $C^\bullet_{CE}(\el,M)$ (\cite{CE}), of the Nijenhuis and Richardson bracket (\cite{NR}), and introduce the subcomplex $C^\bullet_>$ we will need. 
\begin{df}\label{CE}
Let $\el$ be a Lie algebra and $M$ a $\el$-module. Their Chevalley-Eilenberg complex  is the graded space $C(\el,M)$ of antisymmetric linear maps from $\el$ to $M$ :  $$C^n(\el,M):= L(\wedge^n \el,M),$$
with differential 

\begin{eqnarray*}
\delta c(l_1,\dots,l_n) & := & \sum_s (-1)^{s} l_s\cdot c(l_1,\dots, \widehat{l_s}, \dots,l_n)\\
    & + & \sum_{i<j} (-1)^{i+j+1} c([l_i,l_j], l_1,\dots, \widehat{l_i}, \dots,\widehat{l_j}, \dots,l_n),
\end{eqnarray*}
where $\widehat{\ }$ indicates the omission of the underneath term.\\

Moreover $C^{\bullet+1}_{CE}(\el,\el)$ comes equipped with the Nijenhuis-Richardson graded Lie bracket
\begin{equation}\label{NRbra}
[P,Q]:=i_PQ-(-1)^{pq}i_QP,
\end{equation}
where $$i_PQ(l_0,\dots,l_{p+q}):=\sum_{\sigma}sign(\sigma)Q(P(l_{\sigma(0)}, \dots,l_{\sigma(q)}),l_{\sigma(q+1)}, \dots, l_{\sigma(p+q)} )$$ with $\sigma$ in $S_{(p,q+1)}$, the set of permutations such that ${\sigma(0)}<\dots<{\sigma(q+1)}$ and ${\sigma(q+2)}< \dots< {\sigma(p+q+1)}$.

\end{df}

\begin{rem}\label{ad}
When $M=\el$  \cite{NR} have remarked that $d=ad_{\rho_\el}$, where $\rho_\el$ is the bracket of $\el$ seen as an element of $C^{1+1}_{CE}(\el,\el)$.
\end{rem}

\begin{lem} \label{L}
Let $\g$ and $\h$ be two Lie algebras. $\h$ is a $\g\oplus\h$-module via the adjoint action of $\h$ on itself. 
One defines the complex $C_>(\g\oplus\h,\h)$ as the subcomplex of $C(\g\oplus\h,\h)$ defined by
$$C_>(\g\oplus\h,\h)\simeq\bigoplus_{(m,n)\in \mathbb{N}^*\times\mathbb{N}} C^{m,n},$$
where the projection on $C^{m,n}$ is given by restricting elements of $C(\g\oplus\h,\h)$ to $\wedge^m\g\otimes\wedge^n\h$.
It is  a sub-diferential graded Lie algebra of $C(\g\oplus\h,\g\oplus\h)$ endowed with the Nijenhuis-Richardson bracket. Its degree 0 part is abelian and consists of ad-nilpotent elements. We will denote by $(L,d,[\ ,\ ])$ this differential graded Lie algebra. 
\end{lem}

\begin{rem}\label{abus}
In the above lemma we implicitely used the canonical isomorphism $$\wedge^\bullet(\g\oplus\h)\simeq \wedge^\bullet\g\otimes\wedge^\bullet\h,$$
We will denote by the same symbol an element and its image by this isomorphism. When the algebras are finite dimensional, one has $$L\simeq C(\h,\h)\otimes C^+(\g,\mathbb{K}),$$
where the first term is a dgLa, while the second is a commutative graded algebra. 
\end{rem}

\begin{proof}
We start by showing that $C_{>}$ is closed under $d$. Our starting point is, thanks to remark \ref{ad}, that 
$d=ad_{\rho_\g}+ad_{\rho_\h}$.
 But one has $$ad_{\rho_\g} (C^{m,n}) \subset C^{m+1,n}.$$
Indeed, if $c\in C^{m,n}$, one has, since $\rho_\g\in\wedge^2\g^*\otimes\g$, $$i_{\rho_\g}c\in C^{m+1,n}$$ and since $\rho_\g $ admits only terms in $\g$ one has clearly $$i_{c}\rho_\g=0\in C^{m+1,n}.$$
Similarly one can show that 
$$ad_{\rho_\h} (C^{m,n}) \subset C^{m,n+1}.$$
Therefore $$d(C_{>})\subset C_{>}.$$

We now prove that $C_{>}$ is closed under $[\ ,\ ]$.
Let $c_1\in C^{m_1,n_1}$ and $c_2\in C^{m_2,n_2}$. One has immediately that $$i_{c_1}c_2\in C^{m_1+m_2,n_1+n_2-1}.$$
Then $$[c_1,c_2]\in C^{m_1+m_2,n_1+n_2-1}.$$ 
In particular, any element is ad-nilpotent.

Clearly $L(\g,\h)=C^{1,1}$ is abelian. From the above we see that it consists of ad-nilpotent elements.
\end{proof}

\begin{lem}\label{MCL}
Let $\g$ and $\h$ be Lie algebras on $G$ and $H$, then $$J^H_{GGG}+J^H_{GGH}+J^H_{GHH}=0 \iff \rho^H_{GG}+\rho^H_{GH}\in MC(L).$$ 
 \end{lem}

\begin{proof}
Let $c=\rho^H_{GG}+\rho^H_{GH}\in MC(L)$. We analyze the equality $$(\delta c+\frac{1}{2}[c,c])(e_1,e_2,e_3)=0,$$

where $e_i:=g_i+h_i\in \g\oplus\h.$

By definition,
\begin{eqnarray*}
\delta c(e_1,e_2,e_3) & = & - \sum_{\sigma \in\circlearrowright}  c (e_{\sigma(1)},[e_{\sigma(2)},e_{\sigma(3)}]_{\g+\h}) + e_{\sigma(1)}\cdot c (e_{\sigma(2)},e_{\sigma(3)})
\end{eqnarray*}
but 
\begin{eqnarray*}
  c (e_{\sigma(1)},[e_{\sigma(2)},e_{\sigma(3)}]_{\g+\h})  & = & \underset{\heartsuit_1}{\underbrace{\rho_{GG}^H(g_{\sigma(1)},[g_{\sigma(2)},g_{\sigma(3)}]_\g)}}+\underset{\Box_1}{\underbrace{\rho_{GH}^H(g_{\sigma(1)},[h_{\sigma(2)},h_{\sigma(3)}]_\h)}}\\
  & - & \underset{\Delta_1}{\underbrace{\rho_{GH}^H([g_{\sigma(2)},g_{\sigma(3)}]_\g, h_{\sigma(1)})}},
  \end{eqnarray*}

 \begin{eqnarray*}
 e_{\sigma(1)}\cdot c (e_{\sigma(2)},e_{\sigma(3)}) & = &  \underset{\Delta_2}{\underbrace{[h_{\sigma(1)},\rho_{GG}^H(g_{\sigma(2)},g_{\sigma(3)})]_\h}} +\underset{\Box_2}{\underbrace{[h_{\sigma(1)},\rho_{GH}^H(g_{\sigma(2)},h_{\sigma(3)})]_\h}}\\
 & - & \underset{\Box_3}{\underbrace{[h_{\sigma(1)},\rho_{GH}^H(g_{\sigma(3)},h_{\sigma(2)})]_\h}} ,
\end{eqnarray*}

 and

 \begin{eqnarray*}
\frac{1}{2}[c,c](e_1,e_2,e_3) & = &  - \sum_{\sigma \in\circlearrowright} c(e_{\sigma(1)}, c (e_{\sigma(2)},e_{\sigma(3)}))
\end{eqnarray*}

with 

 \begin{eqnarray*}
c(e_{\sigma(1)}, c (e_{\sigma(2)},e_{\sigma(3)})) & = &  \underset{\heartsuit_2}{\underbrace{\rho^H_{GH}(g_{\sigma(1)},\rho^H_{GG}(g_{\sigma(2)},g_{\sigma(3)}))}}\\
& + &  \underset{\Delta_3}{\underbrace{\rho^H_{GH}(g_{\sigma(1)},\rho^H_{GH}(g_{\sigma(2)},h_{\sigma(3)}))}}
\end{eqnarray*}

Regrouping the terms marked with $\Delta$ (resp $\heartsuit$, $\Box$) gives the expression of $J_{GGH}^H$ (resp $J_{GGG}^H$, $J_{GHH}^H$) which was computed in the proof of lemma \ref{Jac}. Note that $\Delta_3 =\Delta_3^{'}+\Delta_3^{''}$ thanks to remark \ref{abus}.

\end{proof}

\begin{cly}\label{Z2}
 $$Z^2_{nab}(\g,\h)\simeq MC(L).$$
\end{cly}

\begin{proof}
By proposition \ref{class}, we know that an element in $Z^2_{nab}(\g,\h)$ is the same thing as an extension $\e$ of $\g$ by $\h$. But remark \ref{form} tells that it boils down to the fact that the Jacobiator of the bracket $\rho :=  \rho_{GG}^G+\rho_{GG}^H+\rho_{GH}^H+\rho_{HH}^H $  vanishes. 

Since $\g$ and $\h$ are supposed to be Lie algebras, proposition \ref{Jac} tells that this is equivalent to the vanishing of $J^H_{GGG}+J^H_{GGH}+J^H_{GHH}$, which in turn means (lemma \ref{MCL}) that  $\rho_{GG}^H+\rho_{GH}^H \in MC(L).$
\end{proof}

%%%%%%%%%%%%%%%%%%%%%%%%%%%%%%%%%%%%%%%%%%%%%%%%%%%%%%%%%%%%%%%%%%%%%%%%%%%%%%%%%%%%%%%%%%%%%%%%%%%%%%%%%%%%%%%%%%%%%%%%%%%%%%

\subsection{Non abelian cohomology as  Deligne groupoid}\label{gauge}

\begin{prop} \label{equiv}
Equivalence in $Z^2_{nab}(\g,\h)$ coincides with gauge equivalence in $MC(L)$.
\end{prop}

\begin{proof}
By lemma \ref{L}, it is legitimate to consider the Deligne groupoid of L.
We recall, by definition \ref{equiMC} and formula (\ref{gbeta}), that two elements $c$ and $c'$ in $MC(L)$ are equivalent if there exists $\beta$ in $\g^*\otimes\h$ such that $$c'=e^{ad_\beta}c+g_\beta.$$  
with $$g_\beta:= - \sum_{n\in \mathbb{N}} \frac{1}{(n+1)!}{ad_\beta}^{n}d\beta.$$

Let us consider $c:=\psi+\chi$, and denote $e_i=g_i+h_i\in \g\oplus\h$.
\begin{eqnarray*}
e^{ad_\beta}(\psi+\chi)(e_1,e_2) & = & (\psi+\chi+[\beta,\psi+\chi] + \underset{0}{\underbrace{ \frac{1}{2} [\beta,[\beta,\psi+\chi]] + \dots}}) (e_1,e_2)\\
& = &\psi(e_1,e_2)+\chi(e_1,e_2)+ \psi(\beta(g_1),g_2)+ \psi(g_1,\beta(g_2))\\
& = & \psi(g_1,h_2)-\psi(g_2,h_1)  +\chi(g_1,g_2)-\psi_{g_2}(\beta(g_1))+\psi_{g_1}(\beta(g_2)),
\end{eqnarray*}
We now compute $g_\beta$. One has
\begin{eqnarray*}
-d\beta(e_1,e_2) & = & -[\rho_\g+\rho_h,\beta](e_1,e_2)\\
& = & -\beta([g_1,g_2]_\g)+[\beta(g_1),h_2]_\h +[h_1,\beta(g_2)]_\h
\end{eqnarray*}
and 
\begin{eqnarray*}
-[\beta,d\beta](e_1,e_2) & = &+2[\beta(g_1),\beta(g_2)]_\h.
\end{eqnarray*}
More generally, for $n>1$
\begin{eqnarray*}
(ad_\beta)^n(d\beta) & = & 0.
\end{eqnarray*}
Therefore, one has
$$ g_\beta=-\beta([g_1,g_2]_\g)+[\beta(g_1),h_2]_\h +[h_1,\beta(g_2)]_\h +[\beta(g_1),\beta(g_2)]_\h.$$
Combining these results gives
\begin{eqnarray*}
c'(e_1,e_2) & = & c(e_1,e_2)+ad^\h_{\beta(g_1)}(h_2) -ad^\h_{\beta(g_2)}(h_1)+ \psi_{g_1}(\beta(g_2))-\psi_{g_2}(\beta(g_1))\\
& + & [\beta(g_1),\beta(g_2)]_\h-\beta([g_1,g_2]_\g).
\end{eqnarray*}
\\

On the other hand, seen as elements in $Z^2_{nab}(\g,\h)$, the two cocycles $c$ and $c'$  are equivalent if equations (\ref{eq1}) and (\ref{eq2}) are satisfied, i.e., using remark \ref{abus}, if  \begin{eqnarray*}
c'(e_1,e_2) & = & c(e_1,e_2)+ad^\h_{\beta(g_1)}(h_2) -ad^\h_{\beta(g_2)}(h_1)+ \psi_{g_1}(\beta(g_2))-\psi_{g_2}(\beta(g_1))\\
& + & [\beta(g_1),\beta(g_2)]_\h-\beta([g_1,g_2]_\g).
\end{eqnarray*}

hence the proposition is proved.

\end{proof}

%%%%%%%%%%%%%%%%%%%%%%%%%%%%%%%%%%%%%%%%%%%%%%%%%%%%%%%%%%%%%%%%%%%%%%%%%%%%%%%%%%%%%%%%%%%%%%%%%%%%%%%%%%%%%%%%%%%%%%%%%%%%%%

\section{Abelian cohomology as  tangent complex}\label{tangent}

We now want to explain how the classification of abelian extensions of Lie algebras fits into the picture given by Theorem \ref{princip}. 

We start by recalling the notion of abelian extension. 

\begin{df}
Let $\g$ be a Lie algebra and let $H$ be a $g$-module. An abelian extension $\e$ of $\g$ by $H$ is a short exact sequence of the form
$$0\to H\to \e \to \g \to 0.$$
\end{df}

Such extension are also considered modulo equivalence :

\begin{df}
Let $\e$ and $\e'$ be two extensions of $\g$ by $H$. They are said to be equivalent if there exists a commutative diagram 

  $$\xymatrix{
      0 \ar[r]  & H \ar[r] \ar@{=}[d]& \e \ar[r] \ar[d]_\varphi&  \g  \ar[r] \ar@{=}[d]& 0  \\
   0 \ar[r] & H \ar[r] & {\e^{'}} \ar[r] &  \g  \ar[r] & 0. }$$
  
\end{df}

Analogously to proposition \ref{class} of the non-abelian case, one has, in terms of the Chevalley-Eilenberg cohomology recalled in definition   \ref{CE} :

\begin{prop}
Abelian extensions of the Lie algebra $\g$ by the $\g$-module $H$ are classified, modulo equivalence, by $H^2(\g,H)$.
\end{prop}

\begin{proof}
This resust is standard \cite{CE} and can be seen as a special case of proposition \ref{class}. One simply needs to read again section \ref{artin}  keeping in mind that now $\psi :=\rho_{GH}^H$ is fixed, defines a $\g$-module structure on $H$ and that $[\ , \ ]_\h=0$. In particular, with these constraints equation $(\ref{mod})$ is  automatically satisfied, while equation $(\ref{cocy})$ means that $\chi$ is a Chevalley-Eilenberg 2-cocycle.  Equation $(\ref{eq1})$ means that $\psi$ is constant while equation $(5)$ now means that $\chi$ and $\chi'$ differ by a Chevalley-Eilenberg coboundary.
\end{proof}

We recall the notion of twisting $\g_\alpha$ of  a dgla $\g$ by a Maurer-Cartan element $\alpha$, sometimes also called tangent complex at $\alpha$.

\begin{df}
Let $\g=(V,[\, \ ], d)$ be a differential graded Lie algebra and let $\alpha\in MC(\g)$ be one of its Maurer-Cartan elements. The tangent complex of $\g$ at $\alpha$ is the differential graded Lie algebra  $\g_\alpha=(V,[\, \ ]_\alpha, d_\alpha)$ where 
$$[u,v]_\alpha=[u,v]$$ 
and 
\begin{equation}\label{dalpha}
d_\alpha(u)=d(u)+[\alpha,u],
\end{equation}
 for $u$ and $v$ in $V$.
\end{df}

\begin{prop}\label{tangentcomplex} Let $\g$ be a Lie algebra and H be a $\g$-module. We denote by $\alpha=\rho_{GH}^H$ the module structure $\rho_{GH}^H: \g\times H\mapsto H,$ or more precisely (remark \ref{abus}) its anti-symetrization. 
With the notations of lemma \ref{L} and definition \ref{CE}, $(C(\g,H),\delta)$ is a sub dgLa of $L_\alpha$.

\end{prop}

\begin{proof}
It is clear that $C(\g,H)$ is an abelian subalgebra of $(L,[\ ,\ ])$. It suffices now to identify  $d_\alpha$ with the Chevalley-Eilenberg coboundary operator $\delta$. One can remark that since $\h:=H$ is abelian, $d$ takes the form 
$$d c(g_1,\dots,g_n) = \sum_{i<j} (-1)^{i+j+1} c([g_i,g_j]_\g, g_1,\dots, \widehat{g_i}, \dots,\widehat{g_j}, \dots,g_n).$$
On the other hand, by definition of the Nijenhuis-Richardson bracket (equation (\ref{NRbra})),
$$[\alpha,c]= \underset{0}{\underbrace{i_\alpha(c)}}-(-1)^{\vert \alpha\vert \vert c\vert}i_c(\alpha),$$
so in particular, since $\vert c\vert=n-2$ and $\vert \alpha \vert=1,$
\begin{eqnarray*}
[\alpha,c](g_1,\dots,g_n) & = & -(-1)^{n-2}\sum_{\sigma\in S_{(1,n-2)}}sign(\sigma)\alpha(c(g_{\sigma(0)}, \dots,g_{\sigma(n-2)}),g_{\sigma(n-1)} )\\
 & = &  -(-1)^{n-3}\sum_{i=1}^n (-1)^{n-i}\alpha (g_i, c(g_1,\dots, \widehat{g_i}, \dots,g_n))\\
 & = &  \sum_{i=1}^n (-1)^{-i}g_i\cdot c(g_1,\dots, \widehat{g_i}, \dots,g_n).
\end{eqnarray*}
Therefore, by definition of $d_\alpha$ (equation (\ref{dalpha})), one has
$$d_\alpha=\delta.$$
\end{proof}

\noindent\textbf{Acknowledgements:}  {We want to thank Boris Shoikhet for having pointed us a mistake, Eugen Lerman and Alan Weinstein for a remark on terminology. We also wish to thank the the anonymous referee for his suggestions of improvement of the text. Support from  Marie Curie Grant IOF-hqsmcf-274032 and SNF 200020\_149150/1.}

%%%%%%%%%%%%%%%%%%%%%%%%%%%%%%%%%%%%%%%%%%%%%%%%%%%%%%%%%%%%%%%%%%%%%%%%%%%%%%%%%%%%%%%%%%%%%%%%%%%%

\bibliographystyle{habbrv}
\bibliography{ghostbib}

\end{document}